\documentclass[journal]{IEEEtran}
\usepackage[utf8]{inputenc}
\usepackage{cite}
\usepackage{graphicx}
\usepackage{url}
\usepackage{lipsum}
\usepackage{color}
\usepackage{xcolor} 
\usepackage{amsmath}
\usepackage[english]{babel}
\usepackage{amsthm}

\theoremstyle{definition}

\usepackage{mathtools}
\usepackage{array}
\usepackage{multirow}
\usepackage{makecell}
\usepackage{float}
\usepackage{setspace}
\usepackage{supertabular,booktabs}
\usepackage{amssymb}
\usepackage{algorithm}
\usepackage{physics}
\usepackage[shortlabels]{enumitem}
\usepackage{algpseudocode}
\usepackage{comment}
\newboolean{showcomments}
\setboolean{showcomments}{true}

\usepackage{amsfonts}
\usepackage{multicol}
\usepackage[top=1.2 cm, bottom=1.2 cm, left=1.35 cm, right=1.35 cm]{geometry}

\allowdisplaybreaks[4]
\usepackage{authblk}
\usepackage{nomencl}
\usepackage[symbol]{footmisc}

\usepackage{subcaption}
\newcommand{\changed}[1]{\textcolor{black}{#1}}

\newcommand{\julio}[1]{\ifthenelse{\boolean{showcomments}}
        { \textcolor{brown}{(JB:  #1)}}{}}
\newcommand{\brian}[1]{\ifthenelse{\boolean{showcomments}}
        { \textcolor{purple}{(BL:  #1)}}{}}
\makenomenclature
\usepackage[colorlinks=true,linkcolor=black,anchorcolor=black,citecolor=black,urlcolor=black]{hyperref}

\begin{document}

\title{\LARGE Computation of  Robust Dynamic Operating Envelopes Based on Non-convex OPF for Unbalanced Distribution Networks}

\author{Bin Liu,~\IEEEmembership{Member,~IEEE} and Julio H. Braslavsky,~\IEEEmembership{Senior Member,~IEEE}
	\thanks{Bin Liu (\textit{corresponding author}) was with Energy Centre, Commonwealth Scientific and Industrial Research Organisation (CSIRO), Newcastle 2300, Australia and is now with the Network Planning Division, Transgrid, Sydney 2000, Australia. Julio, H. Braslavsky is with Energy Centre, CSIRO, Newcastle 2300, Australia. (Email: eeliubin@ieee.org, julio.braslavsky@csiro.au)}}

\markboth{Journal of \LaTeX\ Class Files,~Vol.~22, No.~22, December~2022}%
{Shell \MakeLowercase{\textit{et al.}}: A Sample Article Using IEEEtran.cls for IEEE Journals}

\maketitle

\begin{abstract}
  Robust dynamic operating envelopes (RDOEs) solve the problem of secure allocation of latent network capacity to flexible distributed energy resources (DER) in unbalanced distribution networks. As the computational complexity of RDOEs is much higher than that of \changed{dynamic operating envelopes (DOEs)}, which disregard uncertainties in network parameters and DER capacity utilisation, existing approaches to computing RDOEs have relied on linearised unbalanced three-phase optimal power flow (UTOPF) models to numerate the network feasible region approximately. The use of linearised models, however, risks producing RDOEs that undermine network integrity due to inherent errors in the approximation. This letter presents a practical sensitivity-filtering technique to simplify RDOE numerical computation based on non-convex UTOPF formulations. The accuracy and efficiency of the proposed approach are demonstrated on RDOE allocation with various fairness metrics by testing on representative Australian distribution networks.
\end{abstract}


\section{Introduction}
Dynamic operating envelopes (DOEs) are a key enabler in emerging power system architectures for flexible DER integration \cite{DEIP2022,Liu2021-doe}. It specifies ranges for customers' DER power imports and exports that can be allocated within secure network operations. \changed{Robust DOEs (RDOEs)} extend DOEs to account for unavoidable uncertainties in DOE computation arising from unpredictable utilisation of DOE capacity \cite{liu2022robust,BL_superESD,linearRDOE} and errors in estimated network impedances and forecasts of unmanaged loads \cite{liu2022robust_02,Yi2022}.

An RDOE formulation hedged against uncertainties in customers' utilisation of their allocated capacities was introduced in \cite{liu2022robust} along with a computation methodology. The optimality of the approach was subsequently improved in a superellipsoid formulation in \cite{BL_superESD}.  Improvements in computational efficiency through linear programming techniques were proposed in \cite{linearRDOE}.
These RDOE approaches rely on \emph{linearised} unbalanced three-phase optimal power flow (UTOPF) formulations to make the RDOE computation tractable. This, however, comes at the expense of introducing potential risks to network integrity due to errors inherent to the linear approximation (for illustration, note the differences between network feasible regions (FRs) obtained from linear and non-convex UTOPF formulations in Fig.~\ref{fig_ncvxRDOEgeomaxS_5_bus_export_withQfalse}). 

\changed{
The focus of this letter is to eliminate these potential risks by dropping linear approximations altogether. Instead, we introduce a technique to enable numerically efficient computation of RDOEs based on full nonlinear, non-convex UTOPF formulations. In summary, the main contributions are:}
\begin{enumerate}[leftmargin=*]
\item  \changed{A numerical sensitivity-filtering technique to identify worst-case DOE utilisation scenarios and simultaneously determine critical nodes where network constraints from customers' exports/imports are most likely binding.}

\item \changed{An approach similar to Stochastic Optimisation (SO) to efficiently compute RDOEs based on a non-convex UTOPF formulation ---  the first in this research area --- building upon the previously identified worst-case DOE utilisation scenarios.} 
\end{enumerate}

The approach's numerical efficiency is tested in a case study of optimisation objective options to balance maximum aggregate RDOE capacity versus fairness of RDOE capacity allocation,  demonstrated on a set of representative low-voltage Australian distribution networks. A highlight observation from the study is that a fairer allocation strategy is not necessarily overly conservative with respect to the (typically inherently unfair) allocation strategy that maximises aggregated RDOE capacity.

\section{Non-convex UTOPF-based RDOE}
\subsection{Mathematical formulation}
As discussed in \cite{liu2022robust}, calculating RDOEs is equivalent to seeking the maximum hyperrectangle within the network FR of all DOE customers, where each edge of the hyperrectangle represents the range of active power imports and exports of each DOE customer. We in this letter propose to calculate the RDOE while constraining all vertices of the hyperrectangle within the FR, where each vertex maps to a worst-case DOE utilisation scenario. Using a common non-convex UTOPF \changed{current-voltage (IV)} formulation (see, e.g., Equation~(1) in \cite{liu2022robust}) the RDOE calculation can now be compactly written as
\begin{subequations}\label{rdoeC}
	\begin{eqnarray}
     \max_{p^+,p^-,q}{s(p^+ - p^-)}\\
     \label{rdoeC_01}
          s.t.~~Ap^m+Bq+C\, t(v^m,l^m)=b~\forall m
          \\
            Dv^m+El^m=d,~
            g(v^m)\le f~\forall m\\
			\label{rdoeC_02}
            p^+_k\ge p^-_k,~
            \delta_k^+p^+_k+\delta_k^-p^-_k=0~\forall k
	\end{eqnarray}
\end{subequations}
where $s(p^+-p^-)$ is a suitable objective function of the customers' import-export potential capacity range $p^+ - p^-$, denoted by $\Delta p=p^+-p^-$ in the sequel, where $p^+$ and $p^-$ are the vectors of all DOE customers' upper and lower capacity limits. 
\color{black} The vector $q$ represents a vector consisting of all controllable reactive powers from customers.
The worst-case DOE utilisation scenarios are represented by $p^m~(\forall m)$, which denotes the vertices of a \emph{hyperrectangle} spanned by the scenarios indexed by the superscript $m$. For example, if $p^m$ is a scenario where $K$ DOE customers are exporting at their allocated capacity limits, we have $p^m=[p^-_1,\cdots,p^-_k,\cdots,p^-_K]^T$.

\color{black}
State variable vectors for line currents and nodal voltages are denoted by $l$ and $v$, and $A,B,C,b,D,E,d$ and $f$ are constant parameters with appropriate dimensions. The compact expression \eqref{rdoeC_01} represents the usual relationships between nodal active and reactive powers and products $t(v,l)$ of currents and voltages in current-voltage (IV) OPF formulations. \changed{The equality $Dv+El=d$ represents the voltage drop along each line, and the inequality $g(v)\leq f$ models upper and lower voltage magnitude limits. The discrete variable pair $(\delta_k^+,\delta_k^-)\in \{(1,0),(0,1),(1,1)\}$ indicates if the customer $k$ is exporting power, importing power, or in unknown operational status, and such settings ensure the optimisation results are in line with what is required by a customer. For example, $(\delta_k^+,\delta_k^-)=(1,0)$ forces $p^+_k=0$ while $p^-_k$ is to be optimised, and the allocated RDOEs allow the customer to freely vary its exported power between 0 kW and $|p^-_k|$ kW without violating operational constraints.
Similarly, when $(\delta_k^+,\delta_k^-)=(1,1)$, we have $p_k^+=-p_k^-$, implying the allocated export and import limits for this customer are identical to each other, and it can export and import power up to $|p^-_k|$ kW and $p^+_k$, respectively. Other flexibility can also be achieved by modifying the constraint, and one case is that if the maximum export and import powers are not required to be identical but are limited to specific numbers, say $a$ kW and $b$ kW respectively, \eqref{rdoeC_02} can be updated as $a\ge p^+_k\ge p^-_k\ge -b$.}
For customers with controllable reactive power (DOE customers), $q$ is an optimisation variable, while active and reactive powers for other customers are treated as parameters of the optimisation problem. 

\color{black}
It should be noted that a widely used SO formulation usually models the uncertainties of one or more parameters by a set of scenarios (with the probability associated with each of them), which are randomly sampled before solving the optimisation problem. By contrast, scenarios need to be generated for $p$, which is also a variable, in the proposed approach. However, since $p^m=[p^m_1,\cdots,p^m_k,\cdots,p^m_K]^T$ can be further expressed as $\mathcal{D}(\mathcal{\bar H}^+_{m,:})p^++\mathcal{D}(\mathcal{\bar H}^-_{m,:})p^-$ with $\mathcal{D}$ being a diagonalising operator, and $\mathcal{\bar H}^+_{m,:}$ and $\mathcal{\bar H}^-_{m,:}$ can be determined beforehand\footnote{How $\mathcal{\bar H}^+_{m,:}$ and $\mathcal{\bar H}^-_{m,:}$ are determined will be discussed in the next subsection.}, the proposed approach still employs SO techniques. In this regard, the approach is also similar to the \emph{column-and-constraint generation} method in solving two-stage robust optimisation problems in \cite{ZENG2013457}, which, however, differs in that a sensitivity-filtering approach is employed instead of solving a \emph{max-min} optimisation problem iteratively, to identify the worst scenarios due to the different natures of the problems investigated.

\color{black}
It is also noteworthy that the proposed approach assumes that the intersection of the non-convex FR for DOE customers and the hyperrectangle defining the RDOE is a subset of the former one. A rigorous validation of this assumption falls beyond the scope of the present letter, where the equivalence is numerically tested via post-calculation evaluation in the case studies.

RDOE computation based on the full non-convex UTOPF \eqref{rdoeC} becomes rapidly intractable when the number of DOE customers $K$ is large. Even considering only DOE utilisation scenarios defined by combinations of capacity range limits, i.e. $p^-$ and $p^+$, for each DOE customer, produces a total number of $2^K$ DOE utilisation scenarios for \eqref{rdoeC}. However, the number of non-redundant scenarios that need to be considered can be much smaller. For example, for a distribution network with two DOE customers, of all the $2^2=4$ scenarios, i.e., $(p^+_1,p^+_2)$, $(p^+_1,p^-_2)$, $(p^-_1,p^+_2)$ and $(p^-_1,p^-_2)$, considering only $(p^+_1,p^+_2)$ and $(p^-_1,p^-_2)$ is sufficient if the two customers are at the same bus and the same phase. This motivates us to identify the worst-case DOE utilisation scenarios before solving the optimisation problem, as discussed next.

\subsection{Identifying worst-case DOE utilisation scenarios} 
 Worst-case DOE utilisation scenarios are identified based on nodal/phase voltage sensitivity to DOE customers' demands, which is numerically computed via perturbations introduced in each DOE customer's demand in this letter.
 \color{black}Specifically, from a known network operating condition determined by the demands of all customers, say $(p_0,q_0)$, voltage magnitudes of all bus/phase pairs can be recorded after running unbalanced three-phase power flow (UTPF) via tools such as \texttt{PowerModelsDistribution.jl} \cite{pmd_ref}. Then, the active demand at customer $k$ will be increased by $\Delta w$, a positive and small perturbation, while the demands of all the other customers are kept unchanged, and voltage magnitudes of all bus/phase pairs are updated after running UTPF again. Accordingly, the sensitivity of voltage magnitude at phase $\phi$ of bus $k$ to $p_i$ at the current operating point can be calculated as 
\begin{eqnarray}
    \beta_{i,k}^{\phi}=\frac{\partial u_{i,\phi}}{\partial p_k}|_{(p_0,q_0)}
    \approx\frac{u_{i,\phi}|_{p=p_{-k}}-u_{i,\phi}|_{p=p_0}}{\Delta w}
\end{eqnarray}
where $p_{-k}$ is a vector with all entries the same as $p_0$ except that its $k$-th entry is increased by $\Delta w$; $\phi$ indexes the node phase with $\phi\in\{a,b,c\}$, and $u_{i,\phi}$ is the voltage magnitude of phase $\phi$ at bus $i$.

With the calculated sensitivities, the worst-case DOE utilisation scenarios need to be further identified. As $\beta_{i,k}^{\phi}<0$ implies that $u_{i,\phi}$ becomes lower with increasing $p_k$ and vice visa, $p_k=p^+_k$ is likely to result in lowest $u_{i,\phi}$ and $p_k=p^-_k$ to highest $u_{i,\phi}$. Therefore, to ensure $u_{i,\phi}$ is, for example, under its upper limit in scenario $p^m$, the $k$-th entry of $p^m$ will be set as $p^-_k$. The other entries for $p^m$ and entries for other worst-case DOE utilisation scenarios corresponding to other voltage constraints can be determined similarly. 
 
Further, the worst-case DOE utilisation scenarios can be merged if the values of the calculated sensitivities are sufficiently small. For example, if only voltage upper limits are to be considered for two bus/phase pairs, say phase $a$ at bus $1$ and phase $b$ at bus $1$, in a distribution network with two DOE customers, indexed by $1$ and $2$, and we have $\beta^a_{1,1}>0$, $\beta^a_{1,2}>0$ and $\beta^b_{1,1}>0$, $\beta^b_{1,2}\approx 0$, considering one wort-case scenario $p^1=[p^+_1,p^+_2]^T$, instead of considering both $p^1=[p^+_1,p^+_2]^T$ and $p^2=[p^+_1,0]^T$, is sufficient.

The approach to identifying the worst-case DOE utilisation scenarios, considering all voltage constraints and all DOE customers, is generalised below.

\color{black}
With \changed{$\beta_{i,k}^{\phi}~(\forall i,\forall k)$} for all DOE customers, a $I\times K$ dimension matrix  $\mathcal{H}=[\theta_{i,k}^\phi]_{I\times K}^{\phi\in\{a,b,c\}}$ can be constructed, with $I$ denoting the total number of buses. The entry
$\theta_{i,k}^{\phi}$ takes $1$ when $\beta_{i,k}^{\phi}>\varepsilon$, $-1$ when $\beta_{i,k}^{\phi}<-\varepsilon$ and $0$ when $|\beta_{i,k}^{\phi}|\le\varepsilon$, where $\varepsilon$ is a sufficiently small predefined threshold. 
Then, a new matrix $\mathcal{\bar H}=[h_{mk}]_{M\times K}$ can be derived after removing or merging rows in $\mathcal{H}$ as follows: 
\begin{enumerate}[leftmargin=*]
    \item For any two rows, say $\mathcal{H}_{m,:}$ and $\mathcal{H}_{r,:}$, if there are $\mathcal{H}_{m,j}=\mathcal{H}_{r,j}~(\forall j\in\mathcal{J})$ and $|\mathcal{H}_{m,j}-\mathcal{H}_{r,j}|\le 1~(\forall j\notin\mathcal{J})$ 
  for a non-empty index set $\mathcal{J}$, replace $\mathcal{H}_{m,:}$ and $\mathcal{H}_{r,:}$ by a new row vector, say $\hat h$, where $\hat h_{j}=\mathcal{H}_{m,j}~(\forall j\in\mathcal{J})$ and $\hat h_{j}=\mathcal{H}_{m,j}+\mathcal{H}_{r,j}~(\forall j\notin\mathcal{J})$.
  \item Repeat the above step in the newly derived matrix until any two rows cannot be merged further. 
\end{enumerate}

With $\mathcal{\bar H}$, the constraint \eqref{rdoeC_01} can now be replaced by $A\left(\mathcal{D}(\mathcal{\bar H}^+_{m,:})p^++\mathcal{D}(\mathcal{\bar H}^-_{m,:})p^-\right)+Bq+Ct(v^m,l^m)=b$,
where $\mathcal{\bar H}^+_{m,:}=\max\{\mathcal{\bar H}_{m,:},0\}$ and $\mathcal{\bar H}^-_{m,:}=\min\{\mathcal{\bar H}_{m,:},0\}$. Moreover, sensitivity filtering can be done offline under no-load or estimated operational conditions or updated dynamically based on a network's operational conditions. 

\subsection{Efficiency v.s. fairness in allocating RDOEs}
For the objective function in \eqref{rdoeC}, we take the following allocation strategies for various balances between efficiency and fairness \cite{Kelly1997} and test them on the proposed approach to calculating RDOEs. 

\begin{itemize}[leftmargin=*]
\item Maximum Efficiency (\textit{max\_effcy}), aims to maximise the aggregate RDOE capacity, leading to $s(\Delta p)=\sum\nolimits_k{(\Delta p_k)}$.

\item Proportional Fairness (\textit{ppn\_fair}), with $s(\Delta p)=\sum\nolimits_k{\log\Delta p_k}$, ensures any allocation strategy other than the optimal one, say $\Delta p^*$, has $\sum\nolimits_k{(\Delta p_k-\Delta p^*_k)/\Delta p_k^*}\le 0$.

\item $\alpha$-Fairness (\textit{alpha\_fair}), favouring customers with less allocated RDOE compared with \textit{ppn\_fair}, has $s(\Delta p)=\sum\nolimits_k{-(-\log(\gamma\Delta p_k+\eta))^\alpha}$ with $\alpha$ (a large number that is set as 100 in this letter), $\gamma$ (positive) and $\eta$ being properly selected numbers. Particularly, the max-min fairness, usually regarded as the fairest allocation strategy \cite{Kelly1997}, can be achieved when $\alpha$ approaches $+\infty$. In other words, at the optimal allocation strategy, the RDOE allocation for one customer, say $k$, can not be increased without sacrificing the RDOE allocation, no matter how small, of another customer that already has a poorer allocation than Customer~$k$.

\item Inverse Proportional Fairness to max RDOE (\textit{permax\_fair}), specifies a coefficient for each customer depending on its maximum achievable RDOE (maxPRDOE, denoted as $\Delta p^{ind}_k$)\footnote{To get the maxPRDOE for Customer~$k$, we can set the objective function as $s(\Delta p)=\Delta p_k$.}, leading to $s(\Delta p)=\sum\nolimits_k{\Delta p_k/\Delta p^{ind}_k}$.

\end{itemize}

\changed{It should be noted that although the objective function is convex under each allocation strategy, the allocation model \eqref{rdoeC} is non-convex due to the nature of constraints, which implies that a fair allocation strategy cannot always be guaranteed to exist and, if it exists, can be a local one \cite{brehmer2009proportional}. While the proof of the existence of an allocation strategy and its global optimality is challenging and beyond the scope of this letter, we rely on the solver to determine if a fair allocation strategy exists by reporting the feasibility of the optimisation problem. 
In this letter, all optimisation problems are solved via \texttt{Ipopt} on a desktop machine with Intel(R) Core(TM) i9-9900 CPU and 32 GB RAM.}

While the aggregate RDOE capacity under the \textit{max\_effcy} objective function is expected to be the largest, it is likely also the most unfair. By contrast, \textit{alpha\_fair} is expected to provide the least aggregate RDOE capacity to achieve a fairer allocation.

\section{Case study}\label{case_study}
\begin{figure}[htb!]
	\centering\includegraphics[scale=0.95]{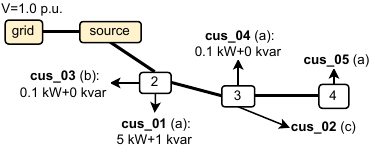}
	\caption{\small The topology and system data for the FiveNetwork (\textit{cus\_02} and \textit{cus\_05} are DOE customers while others are non-DOE customers, and the letter in the bracket indicates the customer's connected phase).}
	\label{fig_5_bus}
\end{figure}

Several distribution networks are studied, including a 5-bus illustrative network (FiveNetwork) and a 33-bus representative Australian network (AusNetwork, whose details are available in \cite{liu2022robust}). The topology and parameters of the FiveNetwork are provided in Fig.\ref{fig_5_bus}. 
The AusNetwork has 87 customers, of which 30 are DOE customers whose reactive powers are also controllable. The active and reactive powers of the remaining non-DOE customers are fixed.
The default limits on active and reactive powers for the FiveNetwork and the AusNetwork are set as [-7 kW, 7 kW] and [-3 kvar, 3 kvar], respectively.

Before presenting the simulation results, the numbers of identified worst-case DOE utilisation scenarios for the two networks and other representative Australian distribution networks from \cite{LVFT_data} are presented in Table \ref{tab_lvft_doe_sce}\footnote{Results for networks with incomplete information are not listed. AusNetwork and Network \texttt{J$^*$} are variations of Network \texttt{J}, where in Network \texttt{J$^*$} all lines are modified to have the same conductor ``3P\_OH\_Code3", and in AusNetwork the distribution transformer connection is changed from delta/wye to wye/wye.}, where bus/voltage sensitivity is approximately estimated by successively increasing customer's demand by 20 kW per phase with $\varepsilon=10^{-5}$. In other words, starting from the non-load network condition, Customer~$k$ is regarded as having an effective impact on $u_{i,\phi}$ only if $|\Delta u_{i,\phi}|$ is higher than $10^{-5}$ p.u. or 0.0023 V for a nominal voltage of 230 V, after setting the demand of Customer~$k$ to 20 kW from 0 kW.

\begin{figure}[htbp!]
	\centering\includegraphics[width=9.5 cm,height=3.7 cm]{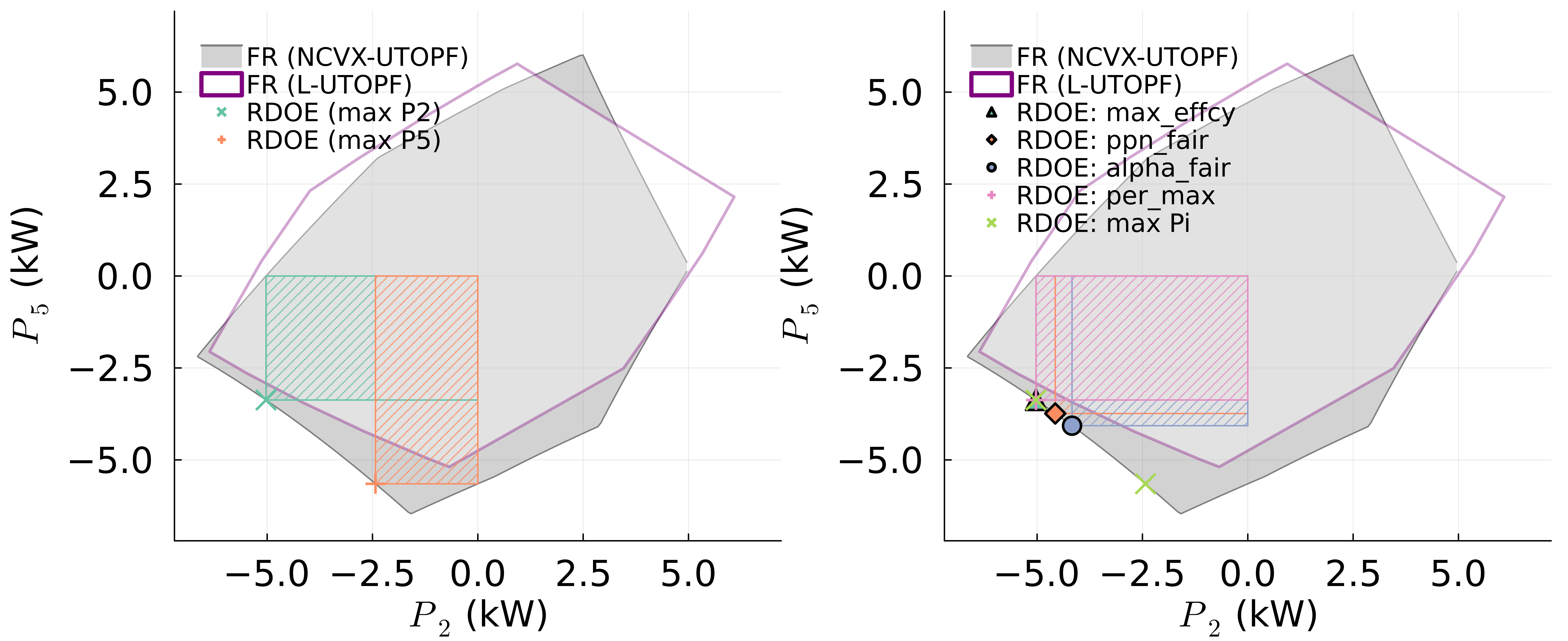}
	\caption{\small Export RDOEs for customers in the FiveNetwork with $q=0~\text{kvar}$ (Left: Maximum Export RDOE achievable by each individual customer; Right: Export RDOEs under various fairness metrics).}
	\label{fig_ncvxRDOEgeomaxS_5_bus_export_withQfalse}
\end{figure}

Results in Table \ref{tab_lvft_doe_sce} show how the total number of worst-case DOE utilisation scenarios (much less than $2^K$)  varies depending on the attributes of the network.
The calculated RDOEs and the FRs for the FiveNework are presented in Fig.\ref{fig_ncvxRDOEgeomaxS_5_bus_export_withQfalse}, where the ``FR (L-UTOPF)" and the ``FR (NCVX-UTOPF)" are calculated according to the approaches in \cite{liu2022robust} and \cite{Riaz2022}, respectively.

\begin{table}[htbp!]\footnotesize\renewcommand\arraystretch{1.00}
	\centering
	\setlength{\tabcolsep}{1.0pt}
	\caption{\small Worst-Case DOE utilisation scenario number (all customers are DOE customers except \texttt{FN} (FiveNework) and \texttt{AN} (AusNetwork), with $\beta^\phi_{i,k}=0~(\forall i)$ if $k$ indicates a non-DOE customer).}
	\begin{tabular}{c|c|c|c|c|c|c|c|c|c|c|c|c|c|c}
		\hline\hline
        Attributes/Network Name & \texttt{FN} & \texttt{AN} & \texttt{G} & \texttt{L} & \texttt{M} & \texttt{N} & \texttt{O} & \texttt{P} & \texttt{R} & \texttt{T} & \texttt{U} & \texttt{V} & \texttt{J} & \texttt{J$^*$} \\\hline
        Total Bus Number & 5 & 32 & 20 & 58 & 82 & 100 & 110 & 158 & 219 & 5 & 11 & 18 & 32 & 32 \\\hline
        Total Cond. Number & 1 & 5 & 4 & 6 & 7 & 10 & 10 & 14 & 12 & 2 & 4 & 5 & 5 & 1 \\\hline
        1-$\phi$ Customer & 5 & 87 & 59 & 36 & 25 & 61 & 82 & 47 & 141 & 2 & 4 & 8 & 87 & 87 \\\hline
        2-$\phi$ Customer & 0 & 0 & 0 & 0 & 0 & 0 & 1 & 2 & 3 & 0 & 0 & 0 & 0 & 0 \\\hline
        3-$\phi$ Customer & 0 & 0 & 0 & 0 & 33 & 2 & 0 & 60 & 0 & 0 & 0 & 0 & 0 & 0 \\\hline
        Worst-Case Scenario Number & 4 & 6 & 20 & 18 & 36 & 6 & 8 & 12 & 12 & 2 & 4 & 14 & 18 & 6 \\\hline
        \hline
	\end{tabular}
	\label{tab_lvft_doe_sce}
\end{table}

As shown in Fig.\ref{fig_ncvxRDOEgeomaxS_5_bus_export_withQfalse}, the use of a linearised UTOPF formulation leads to approximation errors in quantifying the FR that propagates to the RDOE calculations. As expected, the \textit{alpha\_fair} favours Customer~5, which has a smaller RDOE under other allocation strategies.  \begin{table}[htbp!]\footnotesize\renewcommand\arraystretch{1.0}
	\centering
	\setlength{\tabcolsep}{1.0pt}
	\caption{\small Aggregate RDOE capacity and computational times (\texttt{RDOE|Time (seconds)} (``$^*$": $q$ at the optimised value; Otherwise, $q=0$ kvar).}
	\begin{tabular}{c|c|c|c|c|c}
		\hline\hline
		\multirow{2}{*}{Case/Model}  & Linear UTOPF & \multicolumn{4}{c}{NCVX-UTOPF}\\\cline{2-6}
        & ppn\_fair & ppn\_fair & max\_effcy & alpha\_fair &  permax\_fair \\\hline
		FiveNetwork & 7.91$|$0.04& 8.31$|$0.10 & 8.40$|$2.44 & 8.24$|$(0.09) & 8.40$|$(0.10) \\\hline
        FiveNetwork$^*$ & 11.72$|$0.05 & 12.57$|$0.08 & 12.57$|$0.70 & 12.51$|$0.11 & 12.57$|$0.08 \\\hline
		AusNetwork & 138.97$|$4.52 & 137.11$|$0.55& 149.14$|$3.16 & 122.98$|$0.52 & 149.14$|$0.67  \\\hline
        AusNetwork$^*$ & 161.78$|$4.78 & 160.61$|$0.53& 164.66$|$0.57 & 143.87$|$0.80 & 164.66$|$0.55  \\\hline\hline
	\end{tabular}
	\label{tab_case_study}
      \end{table}


Simulation results in Table \ref{tab_case_study} clearly show that, for both the FiveNetwork and the AusNetwork, the \textit{max\_effcy} reports the largest aggregate RDOE capacity, while the \textit{alpha\_fair} leads to the smallest value. It is noteworthy that although Customer~5 is compensated in \textit{permax\_fair}, it reports the same allocated RDOEs as in \textit{max\_effcy}, implying that a compensation mechanism needs to be carefully designed to ensure its effectiveness. The results also clearly demonstrate that RDOEs can be significantly ameliorated via optimising controllable reactive powers.     

\begin{figure}[htb!]
	\centering\includegraphics[scale=0.225]{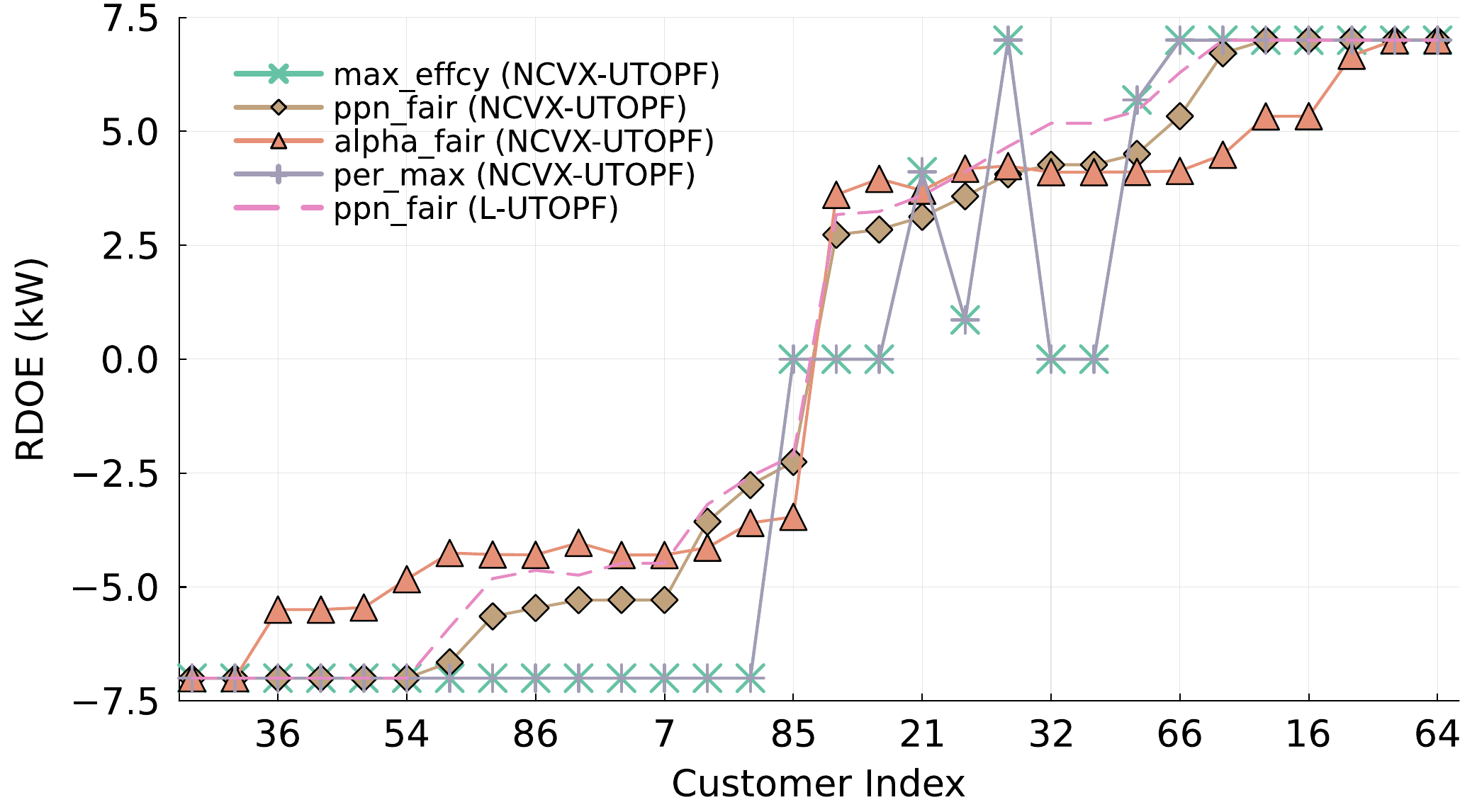}
	\caption{\small RDOEs calculated for customers in the AusNetwork under various fairness metrics with $q=q^\text{opt}$, where negative/positive number indicates the customer is exporting/importing power.}
	\label{fig_ncvxRDOEcurve_aus_J_hybrid_withQtrue}
\end{figure}

The allocation results for DOE customers in the AusNetwork are presented in both Table \ref{tab_case_study} and Fig.\ref{fig_ncvxRDOEcurve_aus_J_hybrid_withQtrue} with controllable reactive powers being optimised. Although the \textit{max\_effcy} (and also \textit{per\_max}) achieves the largest aggregate RDOE capacity, the strategy disallows some customers from exporting or importing powers. By contrast, both the \textit{ppn\_fair} and \textit{alpha\_fair} allocate RDOEs with a stronger favour of \emph{weaker} customers, thus providing a fairer allocation result. \changed{Compared with the approach based on linear UTOPF model, as shown in Fig.\ref{fig_ncvxRDOEcurve_aus_J_hybrid_withQtrue} and similar to the observations in Fig.\ref{fig_ncvxRDOEgeomaxS_5_bus_export_withQfalse}, the proposed approach reports more accurate RDOEs. Specifically, the linear UTOPF-based approach underestimates exporting limits for both the FiveNetwork and the AusNetwork and overestimates export limits in the AusNetwork due to the inevitable errors from linearising the UTOPF model.}

To further verify the effectiveness of the proposed approach, the RDOEs computed for all allocation strategies were evaluated on 30,000 randomly generated operation scenarios for the DOE customers (following the approach suggested in \cite[Section IV.C]{liu2022robust}).
The simulations show that all 30,000 scenarios remain within the network FR, reporting no voltage violations, which demonstrates the approach's robustness.

On computational times, all tested cases are solved within five seconds, as shown in Table \ref{tab_case_study}. The computational efficiency is further tested on the network \texttt{R} in Table \ref{tab_lvft_doe_sce}, where 50 of the 141 customers are set as DOE customers. The sensitivity filtering identifies 6 worst-case DOE utilisation scenarios, and the solver takes 3.72 seconds to report the optimal solution, which demonstrates the efficiency of the proposed approach when applied in a real network on a larger scale. 

\section{Conclusions}
This letter extends the calculation of RDOEs to enable full non-convex UTOPF formulations by exploiting a practical sensitivity-filtering technique to reduce computational complexity. Case studies demonstrate the higher accuracy of the approach as compared to existing methods based on linear UTOPF formulations and show promising computational efficiency in both illustrative and real distribution networks. Trade-offs in efficiency and fairness quantified in the case studies demonstrated that a fairer strategy for robust allocation of network capacity, although leading to a smaller aggregate RDOE capacity, is not necessarily too conservative with respect to the strategy with maximal efficiency. 

\bibliography{REFs_Power_Grid}

\begin{IEEEbiography}[{\includegraphics[width=1.13in, height=1.32in, clip]{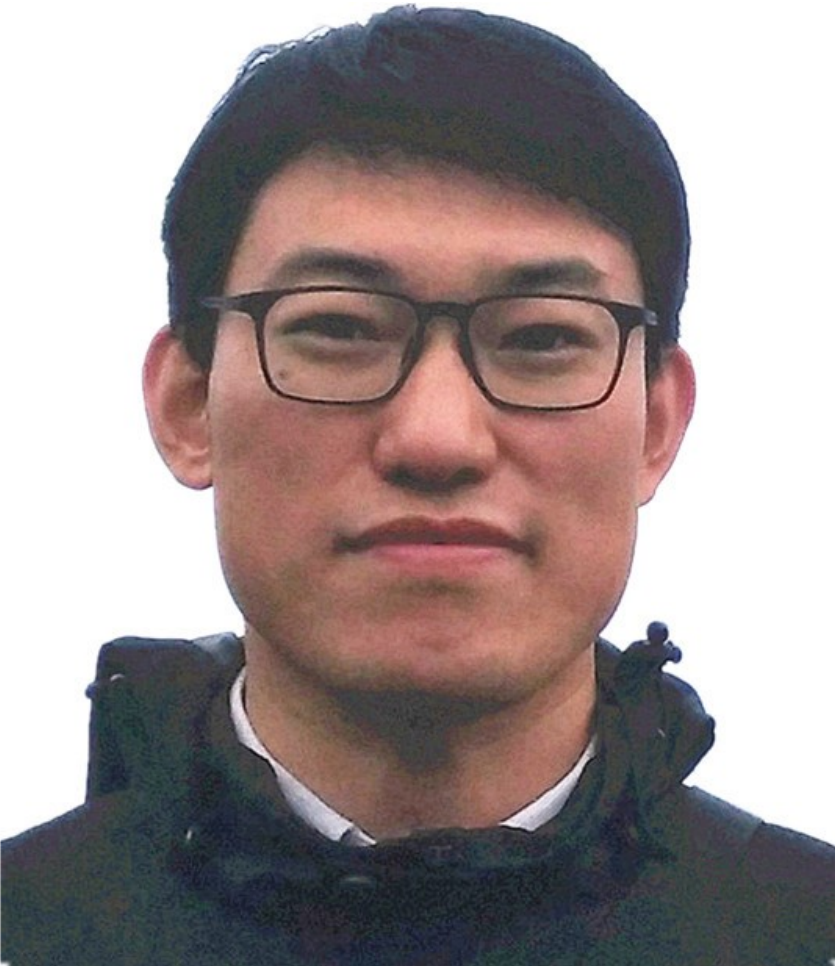}}]{Bin Liu} (Member, IEEE) received his Bachelor's degree from Wuhan University, Wuhan, China, Master's degree from China Electric Power Research Institute, Beijing, China and PhD degree from Tsinghua University, Beijing, China, all in Electrical Engineering, in 2009, 2012 and 2015, respectively. 

He is a Senior Power System Engineer in the Network Planning Division of Transgrid, Sydney, Australia. Before joining Transgrid, he had held research or engineering positions with the Energy Centre of Commonwealth Scientific and Industrial Research Organisation (CSIRO), Newcastle, Australia, The University of New South Wales, Sydney, Australia, the State Grid, Beijing, China, and The Hong Kong Polytechnic University, Hong Kong. His current research interests include power system modelling, analysis \& planning, optimisation theory applications in the power and energy sector, and the integration of renewable energy, including distributed energy resources (DERs).
\end{IEEEbiography}

\begin{IEEEbiography}[{\includegraphics[width=1.13in, height=1.32in, clip]{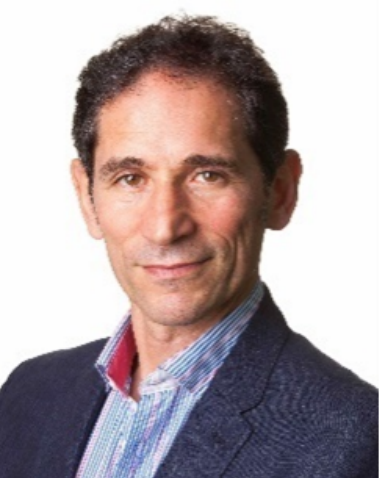}}]{Julio H. Braslavsky}(Senior Member, IEEE) received his PhD in Electrical Engineering from the University of Newcastle NSW, Australia in 1996, and his Electronics Engineer degree from the National University of Rosario, Argentina in 1989. He is a Principal Research Scientist with the Energy Systems Program of the Australian Commonwealth Scientific and Industrial Research Organisation (CSIRO) and an Adjunct Senior Lecturer with The University of Newcastle, NSW, Australia. He has held research appointments with the University of Newcastle, the Argentinian National Research Council (CONICET), the University of California at Santa Barbara, and the Catholic University of Louvain-la-Neuve in Belgium. His current research interests include modelling and control of flexible electric loads and integration of distributed power-electronics-based energy resources in power systems. He is Senior Editor for IEEE Transactions on Control Systems Technology.
\end{IEEEbiography}

\end{document}